\documentclass{amsart}

\usepackage[utf8]{inputenc}
\usepackage[T1]{fontenc}
\usepackage[colorlinks=true]{hyperref}
\usepackage{amssymb}
\usepackage{tikz}
\usepackage{upgreek}

\newcommand\N{\mathbb{N}}
\newcommand\R{\mathbb{R}}
\newcommand{\no}{\ensuremath{\mathbf{No}}}
\newcommand{\on}{\ensuremath{\mathbf{On}}}

\theoremstyle{plain}
\newtheorem{thm}{Theorem}[section]
\newtheorem{prop}[thm]{Proposition}
\newtheorem{cor}[thm]{Corollary}
\theoremstyle{remark}
\newtheorem{rem}[thm]{Remark}
\theoremstyle{definition}
\newtheorem{defn}[thm]{Definition}
\newtheorem{ex}[thm]{Example}

\numberwithin{equation}{section}

\numberwithin{thm}{section}

\title[Surreal numbers with derivation, Hardy fields and transseries]{Surreal numbers with derivation, Hardy fields and transseries: a survey}

\author{Vincenzo Mantova}
\address{University of Leeds}
\email{v.l.mantova@leeds.ac.uk}

\author{Mickaël Matusinski}
\address{University of Bordeaux}
\email{mickael.matusinski@math.u-bordeaux.fr}

\dedicatory{To the memory of Murray Marshall.}

\date{}

\subjclass[2010]{}

\begin{document}

\begin{abstract}
  The present article surveys surreal numbers with an informal
  approach, from their very first definition to their structure of
  universal real closed analytic and exponential field. Then we
  proceed to give an overview of the recent achievements on equipping
  them with a derivation, which is done by proving that surreal
  numbers can be seen as transseries and by finding the `simplest'
  structure of H-field, the abstract version of a Hardy field. All the
  latter notions and their context are also addressed, as well as the
  universality of the resulting structure for surreal numbers.
\end{abstract}

\maketitle

\section{Introduction}
\label{sec:intro}

The theory of surreal numbers initiated by J.H. Conway in
\cite{conway_numb-games}, and popularized by D. Knuth
\cite{knuth:surreal-numbers}, is fascinating and fruitful but not so
well known. However, it has been enhanced continuously and now,
several remarkable achievements have been reached.

Rooted into very fundamental and accessible set theoretic
considerations, these objects nonetheless strike by the richness
together with the universality of their structure. As a class of
numbers -- denote it by \textbf{No} -- they include simultaneously the
real numbers and the ordinal numbers, unified into a common algebraic
structure: a huge real closed field. Moreover, \textbf{No} can be
viewed as a field of (generalized) power series with real
coefficients, being accordingly a -- one can even say ``the'' --
universal domain for all real closed fields. Moreover, after
\cite{gonshor_surreal}, \textbf{No} carries exponential and
logarithmic functions.

These remarkable facts resonated with important results in
differential equations \cite{ecalle:dulac} and model theory
\cite{wil:modelcomp, dmm:exp} concerning tame geometry
(non-oscillating real functions and their formal analogues). Several
similar big real closed fields of formal power series with exponential
and logarithm -- called transseries, or log-exp series, or exp-log
series -- have been developed, but with an important additional item:
a well-behaved derivation. In fact, the purpose of these fields is
precisely to build algebraic structures being closed under resolution
of differential equations in terms of formal analogues of
non-oscillating real functions. This purpose appears to have reached a
milestone recently with the quantifier elimination result in
\cite{adh:mod-theo-transseries}.

Related to this, one of the recent achievements concerning surreal
numbers is the construction in \cite{ber-man:surreal-derivations} of a
well-behaved derivation in \textbf{No}. Moreover, for this derivation
\textbf{No} is closed under integration and even the results in
\cite{adh:mod-theo-transseries} apply: \textbf{No} can be said to be
``the'' universal domain for non-oscillating differentiable real
functions according to \cite{adh:surreal-universal}.

The present survey article is divided into three sections and has two
aims: provide a gentle and intuitive introduction to surreal numbers
and give an overview of these new results about derivation within
their context. In particular, note that we won't address several other
old and new interesting results, e.g.\ the ones about analysis
\emph{on} surreal numbers like
\cite{costin-ehrlich-friedman:integration-surreal} or the ones about
number theory through the notion of integer part like
\cite{ehrlich-kaplan:number-system-II}. To survey these results would
lead us way further and to a much longer article. We will try to stay
guided by the following nice fact. As mathematical objects, surreal
numbers can be viewed in three different manners: as numbers in the
set theoretic sense (cuts between ordered sets: see Definition
\ref{defi:surreal-ind}), as combinatorial objects (ordered binary
sequences: see Definition \ref{defi:surreal-sign}) or as formal
analogues of non-oscillating real functions (generalized power series:
see Theorem \ref{theo:NO-series-gene}). This ubiquity can explain the
possibility for surreal numbers to be used into different contexts, in
particular as an ordered algebraic structure from the model theoretic
and tame geometric point of view. However, we will be more interested
in the third point of view and will provide in the second section a
survey on the corresponding functional and formal objects (Hardy
fields, transseries, H-fields). The last section will be devoted to
survey the results in
\cite{ber-man:surreal-derivations,adh:surreal-universal} and their
corollaries.

The authors would like to thank A.\ Gehret, L.\ van den Dries and an
anonymous referee for the useful comments and corrections.

\section{Surreal numbers}
\label{sec:surreal}

\subsection{Surreal numbers consist of ``all numbers great and
  small''.}
The construction of surreal numbers goes back to
\cite{conway_numb-games}. One of the most striking ideas of
J.H. Conway is to merge two fundamental constructions of numbers:
\emph{Dedekind's construction of real numbers} in terms of cuts in the
set of rational numbers, and \emph{von Neumann's construction of
  ordinal numbers} by transfinite induction in terms of set
membership. All that is needed to get started with surreal numbers is
the basis of set theory (to be precise: the axiomatic ``NEG with
Global Choice'', which extends conservatively ``ZFC'', see
e.g.\ \cite{fitting-smullyan:set-theory-NBG}), the following two
definitions and nothing more than the empty set $\emptyset$. Let us
paraphrase J.H. Conway:
\begin{defn}[Conway]
  \label{defi:surreal-ind}
  $ $
  \begin{itemize}
  \item[\textbf{(S1)}] If $L$, $R$ are any two sets of numbers, and no
    member of $L$ is $\geq$ any member of $R$, then there is a number
    $\{L\,|\,R\}$. All numbers are constructed in this way.
  \item[\textbf{(S2)}] For two numbers $a=\{L\,|\,R\}$ and
    $b=\{L'\,|\,R'\}$, we say that $a\geq b$ iff no member of $R$ is
    $\leq b$ and no member of $L'$ is $\geq a$. Note that $a\leq b$
    just means $b\geq a$.
  \item[\textbf{(S3)}] For two numbers $a, b$, we say that $a = b$ if
    $a \leq b$ and $b \leq a$.
\end{itemize}
\end{defn}

The very first number has to be $\{\emptyset\, |\, \emptyset\}$ which
is naturally called the number $0$. Then come
$1:=\{0\, |\, \emptyset\}$ and $-1:=\{ \emptyset\, |\, 0\}$. We leave
it to the reader to check that these 3 pairs actually define surreal
numbers, and that they verify: $-1<0<1$. Likewise, one obtains at the
next step for example:
\begin{equation}
  \label{equ:example}
  -2:=\{ \emptyset\, |\, -1 \} < -1 < -\frac{1}{2} :=
  \{ -1\, |\, 0\} < 0 < \frac{1}{2} :=
  \{0\, |\, 1 \} < 1 < 2 := \{ 1\, |\, \emptyset \}.
\end{equation}
Actually, just by considering numbers of the type
$\{L\, |\, \emptyset\}$ and by transfinite induction, we can construct
any well-ordered set. E.g., all the natural numbers are recovered via
the simple induction formula: $n+1:=\{n\, |\, \emptyset\}$, and the
first infinite ordinal is:
$$
\omega:=\{1,2,3,\ldots\, |\, \emptyset\}.
$$
In other words, the proper class of ordinal numbers -- denote it by
\textbf{On} -- is included into surreal numbers, the proper class of
the latter being denoted by $\no$. Likewise, anti-well-ordered sets
are defined by restricting to the numbers $\{\emptyset\, |\, R\}$. The
choice of $1/2$ and $-1/2$ for the preceding surreal numbers, with its
implicit arithmetical content, will be justified thereafter.

\medskip

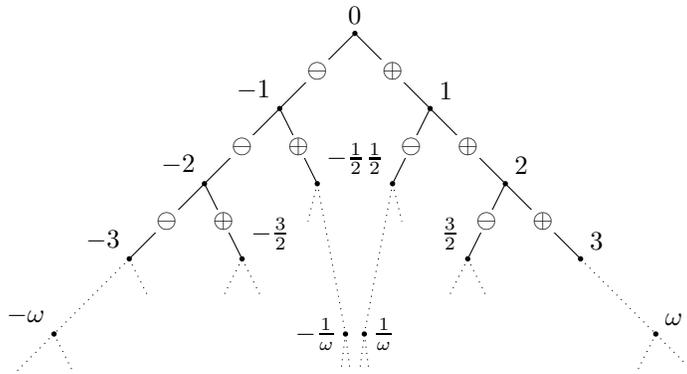
\begin{figure}[b]
  \centering
  \begin{tikzpicture}[radius=1pt,x=1cm,y=1cm]
    \fill (0,0) circle;
    \draw (0,0) node [above] {$0$};

    \fill (-1,-1) circle;
    \draw (0,0) -- (-1,-1)
      node [midway,fill=white,inner sep=.03cm] {$\ominus$};
    \draw (-1,-1) node [above left] {$-1$};
    \fill (1,-1) circle;
    \draw (0,0) -- (1,-1)
      node [midway,fill=white,inner sep=.03cm] {$\oplus$};
    \draw (1,-1) node [above right] {$1$};

    \fill (-2,-2) circle;
    \draw (-1,-1) -- (-2,-2)
      node [midway,fill=white,inner sep=.03cm] {$\ominus$};

    \draw (-2,-2) node [above left] {$-2$};
    \fill (-1/2,-2) circle;
    \draw (-1,-1) -- (-1/2,-2)
      node [midway,fill=white,inner sep=.03cm] {$\oplus$};

    \draw (-1/2,-2) node [above right] {$-\frac{1}{2}$};
    \fill (1/2,-2) circle;
    \draw (1,-1) -- (1/2,-2)
      node [midway,fill=white,inner sep=.03cm] {$\ominus$};

    \draw (1/2,-2) node [above left] {$\frac{1}{2}$};
    \fill (2,-2) circle;
    \draw (1,-1) -- (2,-2)
      node [midway,fill=white,inner sep=.03cm] {$\oplus$};

    \draw (2,-2) node [above right] {$2$};

    \fill (-3,-3) circle;
    \draw (-2,-2) -- (-3,-3)
      node [midway,fill=white,inner sep=.03cm] {$\ominus$};
    \draw [dotted] (-3,-3) -- (-2.75,-3.5);

    \draw (-3,-3) node [above left] {$-3$};
    \fill (-3/2,-3) circle;
    \draw (-2,-2) -- (-3/2,-3)
      node [midway,fill=white,inner sep=.03cm] {$\oplus$};

    \draw (-3/2,-3) node [above right] {$-\frac{3}{2}$};
    \draw [dotted] (-3/2,-3) -- (-7/4,-3.5);
    \draw [dotted] (-3/2,-3) -- (-5/4,-3.5);

    \draw [dotted] (-1/2,-2) -- (-5/8,-2.5);

    \fill (-1/8,-4) circle;
    \draw [dotted] (-1/2,-2) -- (-1/8,-4)
      node [left] {$-\frac{1}{\omega}$};
    \draw [dotted] (-1/8,-4) -- (-3/16,-4.5);
    \draw [dotted] (-1/8,-4) -- (-1/16,-4.5);

    \fill (1/8,-4) circle;
    \draw [dotted] (1/2,-2) -- (1/8,-4)
      node [right] {$\frac{1}{\omega}$};
    \draw [dotted] (1/8,-4) -- (3/16,-4.5);
    \draw [dotted] (1/8,-4) -- (1/16,-4.5);

    \draw [dotted] (1/2,-2) -- (5/8,-2.5);

    \fill (3/2,-3) circle;
    \draw (2,-2) -- (3/2,-3)
      node [midway,fill=white,inner sep=.03cm] {$\ominus$};

    \draw (3/2,-3) node [above left] {$\frac{3}{2}$};
    \fill (3,-3) circle;
    \draw (2,-2) -- (3,-3)
      node [midway,fill=white,inner sep=.03cm] {$\oplus$};
    \draw [dotted] (3/2,-3) -- (7/4,-3.5);
    \draw [dotted] (3/2,-3) -- (5/4,-3.5);

    \draw (3,-3) node [above right] {$3$};

    \draw [dotted] (3,-3) -- (4, -4);
    \fill (4,-4) circle;
    \draw (4,-4) node [above right] {$\omega$};
    \draw [dotted] (4,-4) -- (4.5,-4.5);
    \draw [dotted] (4,-4) -- (3.75,-4.5);

    \draw [dotted] (-3,-3) -- (-4, -4);
    \fill (-4,-4) circle;
    \draw (-4,-4) node [above left] {$-\omega$};
    \draw [dotted] (-4,-4) -- (-4.5,-4.5);
    \draw [dotted] (-4,-4) -- (-3.75,-4.5);
  \end{tikzpicture}
  \caption{The tree of surreal numbers.}
  \label{fig:tree}
\end{figure}

Due to their inductive construction, surreal numbers naturally form a
hierarchy. More precisely, to each surreal number one can assign its
\emph{birthday}: the ordinal number at which it has been constructed
\cite[Theorem 16, p.30]{conway_numb-games}. Moreover, any surreal
number $a$ of birthday $\alpha$ is equal to a cut $\{L\, |\, R\}$
among the numbers with birthday less than $\alpha$. Therefore, for any
ordinal $\beta<\alpha$, $a$ defines a cut among the numbers with
birthday less than $\beta$:
$$
\{b<a,\, \mathrm{ birthday}(b)<\beta\ |\ b>a,\, \mathrm{
  birthday}(b)<\beta\},
$$
the latter cut being itself a number $a_\beta$ of birthday
$\beta$. Following P. Ehrlich's terminology
\cite{ehrlich:all-numb-great-small}, such $a_\beta$ is said to be
\textbf{simpler than} $a$, denoted by $a_\beta<_s a$. The relation
$\leq_s$ is a well-founded partial order relation on $\no$. Conway
comes to this observation after establishing that $\no$ can be viewed
as a \emph{lexicographically ordered binary tree}
\cite[p.11]{conway_numb-games} (see Figure~\ref{fig:tree}).

Then, with the preceding notations, the $a_\beta$'s are exactly the
nodes above $a$, and $a$ itself as a node gives rise at the next stage
to exactly two new surreal numbers, one at its immediate left and
another at its immediate right. Now, to each branch of the tree a sign
can be attributed $\ominus$ if the branch tilts to the left, $\oplus$
if it tilts to the right. This leads us to another representation for
a surreal number: the sequence of signs which expresses the unique
path connecting it to the common root 0. We call it the \textbf{sign
  sequence representation} of a surreal number. Conversely, any sign
sequence of length some ordinal is the sign sequence of a unique
surreal number \cite[Theorem 18]{conway_numb-games}.

\medskip

Gonshor takes the existence of such a representation as a starting
point for the theory of $\no$, assuming that the proper class
\textbf{On} of ordinal numbers is given. In \cite{gonshor_surreal}, he
uses the following equivalent definition for $\no$:
\begin{defn}[Gonshor]
  \label{defi:surreal-sign}
  A \textbf{surreal number} is any function from a given ordinal
  number to the set of two elements: $\{\ominus,\,\oplus\}$, i.e.\
  informally, a sequence of pluses and minuses indexed by some
  ordinal.
\end{defn}
The ordinal $\alpha$ on which the surreal number is defined is called
its \textbf{length} -- set $l(a):=\alpha$ -- which coincides with the
birthday of $a$. The partial ordering on $\no$ called the
\textbf{simplicity ordering} is defined as:
\[
  a\leq_sb \Leftrightarrow \textrm{ the sign sequence of } a \textrm{
    is an initial subsequence of that of } b.
\]

For the total ordering, let $a$ and $b$ be any two surreal numbers
with $l(a)\leq l(b)$. Consider the sign sequence of $a$ completed with
0's so that the two sign sequences have same length. Then consider the
\textbf{lexicographical order} between them, denoted by $\leq$, based
on the following relation:
$$
\ominus<0<\oplus
$$
e.g.\ the inequalities (\ref{equ:example}) can be written as:
$$
\ominus\ominus < \ominus < \ominus\oplus < 0 < \oplus\ominus < \oplus
< \oplus\oplus.
$$

Returning to Conway's construction, the reader might have noticed a
seeming ambiguity. What is for instance the number $\{ -1\, |\, 1\}$?
Or $\{ -2\, |\, 1/2,1,2\}$? In fact, both of them are equal to 0,
i.e.\ $\{\emptyset\, |\, \emptyset\}$. This indicates that the
representation of a number as an actual cut $\{L\,|\,R\}$ is not
unique. The following notion appearing in Gonshor's book
\cite[p.9]{gonshor_surreal} helps to clarify the situation a bit. Take
two pairs $(L,R)$ and $(L',R')$ of \emph{subsets} of $\no$ with $L<R$
and $L'<R'$. $(L',R')$ is said to be \textbf{cofinal} in $(L,R)$ if
for any $(a,b)\in L\times R$, there is $(a',b')\in L'\times R'$, such
that $a\leq a'<b'\leq b$. $(L', R')$ and $(L, R)$ are \textbf{mutually
  cofinal} if $(L', R')$ is cofinal in $(L, R)$ and $(L, R)$ is
cofinal in $(L', R')$. One has that:

\begin{thm}[Cofinality theorem {\cite[Thms.\
    2.6-2.7]{gonshor_surreal}}]
  \label{thm:cofin}
  Let $a = \{L\,|\,R\}$ and $(L', R')$ be cofinal in $(L, R)$. Then:
  \begin{itemize}
  \item if $L' < a < R'$, then $a =\{L'\,|\,R'\}$;
  \item if $(L, R)$ is cofinal in $(L', R')$, namely if $(L, R)$ and
    $(L', R')$ are mutually cofinal, then $a = \{L'\,|\,R'\}$.
  \end{itemize}
\end{thm}

\begin{thm}[Inverse cofinality theorem {\cite[Thms.\
    2.8-2.9]{gonshor_surreal}}]
  Let $a$ be a number, and let
  $L_a:=\{b\in\no\ ;\ b<a\textrm{ and }b<_s a\}$,
  $R_a:=\{b\in\no\ ;\ b>a\textrm{ and }b<_s a\}$. Then
  $a =\{L_a\,|\,R_a\}$. Moreover, if $a =\{L'\,|\,R'\}$, then
  $(L', R') $ is cofinal in $(L_a , R_a)$.
\end{thm}

The simplest surreal number that lies between $-1$ and $1$, or between
$-2$ and $\{1/2,1,2\}$, is indeed 0. For any surreal number $a$, its
representation $\{ L_a\, |\, R_a\}$ for $L_a$ and $R_a$ as above is
called the \textbf{canonical cut} of $a$. By abuse of notation, we
also denote the canonical cut by $a=\{ a^L\ |\ a^R\}$ where $a^L$ and
$a^R$ are general elements of the canonical sets $L_a$ and $R_a$
(e.g.\ $a^L=n$ if $L_a=\mathbb{N}$).

Now the notations $1/2:=\{0\, |\, 1\}$ and $-1/2:=\{1\, |\, 0\}$ can
be partly justified by the following intuitive statement: they are the
simplest numbers between $0$ and $1$, and $-1$ and $0$
respectively. In fact, any surreal number of finite length corresponds
to a \emph{dyadic fraction}. Subsequently, among the surreals of
length $\omega$, we find the rationals non-dyadic and the irrational
real numbers. Note that we also get a priori non-real numbers, like
the simplest infinitely large number $\omega$ (recall that
$\omega:=\{n\, |\, \emptyset\}$), or the simplest infinitesimal one
$\{0\, |\, 1/2^n\}$ to which we assign the symbol $1/\omega$. In the
following section such an assignation will be plainly justified from
an algebraic standpoint.

\subsection{A universal domain for real closed fields}
\label{sec:universal-RCF}
The next tour de force in Conway's construction is not only to
reconstruct numbers as ordered elements but also to recover the
algebraic relations between them, and even to uncover some new ones.

\begin{defn}[\textbf{Field operations}]
  For any surreal numbers $a=\{a^L\,|\, a^R\}$ and
  $b=\{b^L\, |\,b^R\}$, let us define:
  $$
  \begin{array}{lccl}
    \textbf{addition:} & a+b & := & \{ a^L+b,\, a+b^L\, |\, a^R+b,\, a+b^R\}\\
    \textbf{inverse \ element:} & -a & := & \{-a^R\, |\, -a^L\} ; \\
    \textbf{neutral\ element:} & 0 & = & \{\emptyset\ |\ \emptyset\} ; \\
    \textbf{multiplication:} & a \cdot b & := & \{ a^L \cdot b+a \cdot b^L-a^L \cdot b^L,\\
                       &&&\quad a^R \cdot b+a \cdot b^R-a^R \cdot b^R\, |\, \\
                       &&&\quad\ \  a^L \cdot b+a \cdot b^R-a^L \cdot b^R, \\
                       &&&\quad\ \  a^R \cdot b+a \cdot b^L-a^R \cdot b^L\}  \\
    \textbf{neutral\ element:} &1&=&\{ 0\, |\, \emptyset\}=\oplus.
  \end{array}
  $$
\end{defn}
First of all, we need to check that these inductive formulas are
well-defined. For the definitions of the addition and of the inverse
element, one can observe that the elements of the left hand set of the
cut formula are $<$ the elements of the right hand set, for instance
$a+b^L<a^R+b$. For the multiplication, the corresponding necessary
condition is less obvious. It relies on the following type of
observations:
\begin{multline*}
  (a-a^R)(b-b^L)<0<(a-a^L) \cdot (b-b^L) \Leftrightarrow \\
  a^L \cdot b+a \cdot b^L-a^L \cdot b^L<a \cdot b<a^R \cdot b+a \cdot
  b^L-a^R \cdot b^L.
\end{multline*}
Let us illustrate by an example the multiplication formula. Consider
$\omega:=\{n\, |\, \emptyset\}$ and $1/\omega:=\{0\, |\, 1/2^m\}$
(recall that $n$ and $m$ implicitly denote arbitrary positive
integers) and suppose that we have already defined the multiplication
of simpler pairs of numbers. Then:
\begin{align*}
  \omega \cdot \frac{1}{\omega} & = \left\{n\, |\, \emptyset\right\} \cdot \left\{0\, |\,  \frac{1}{2^m}\right\}\\
                                & = \left\{n \cdot \frac{1}{\omega}+0 \cdot \omega-0 \cdot n\, |\, n \cdot \frac{1}{\omega}+\omega \cdot \frac{1}{2^m}-n \cdot \frac{1}{2^m}\right\}\\
                                & = \left\{n \cdot \frac{1}{\omega}\, |\, n \cdot \frac{1}{\omega}+(\omega-n) \cdot \frac{1}{2^m} \right\}\\
                                & = 1 \qquad\qquad \textrm{ since }n \cdot \frac{1}{\omega}<1<(\omega-n) \cdot \frac{1}{2^m}.
\end{align*}

There is no simple formula in terms of sign sequences for the various
operations, except for the additive inverse of a surreal number: the
sign sequence of $-a$ is derived from that of $a$ by changing
$\oplus$'s into $\ominus$'s and $\ominus$'s into $\oplus$'s.

Note that Gonshor proves the so-called \textbf{uniformity properties}
for these operations as well as for the maps that will be defined
later on \cite[Theorems 3.2, 3.5 etc.]{gonshor_surreal}. For example
in the case of the addition, this means that $a+b$ may be obtained by
taking any cuts $a=\{ L\ |\ R\}$ and $b=\{ L'\ |\ R'\}$ in the
respective cofinality classes instead of taking their canonical cuts,
and applying the same formula as above. In other words, \emph{the
  formulas do not depend on the cuts for which $a^L,\ a^R,\ b^L,\ b^R$
  can be taken as general elements}.

\medskip

We already identified real numbers in terms of cuts among the dyadic
fractions:
\begin{defn}
  \label{defi:real-nber}
  $a$ is a real number iff:
  $$
  a=\left\{a - \frac{1}{2^n}\, |\,a+ \frac{1}{2^n}\right\}
  $$
  (where it is understood that $n$ varies in $\mathbb{N}$).
\end{defn}
In terms of sign sequences, Gonshor establishes that real numbers are
those with either finitely many signs or with length $\omega$ and not
ultimately of constant sign. On real numbers, the previously defined
operations restrict well and one thereby obtains a subfield of $\no$
having the least upper bound property, i.e.\ a copy of the field
$\mathbb{R}$.

Addition and multiplication also restrict well to ordinal numbers:
they correspond to the notions of \emph{natural} sum and product, also
called Hessenberg operations. Recall that these operations consist in
taking the \emph{Cantor normal form} of ordinal numbers and applying
addition and multiplication on them as if they were actual polynomials
in an abstract variable $\omega$.

\medskip

Following Conway, a proper class with a field structure is called a
Field. Concerning $\no$ itself, one obtains:

\begin{thm}[{\cite[Ch. 1]{conway_numb-games}, \cite[Theorems 3.3, 3.6,
    3.7 ]{gonshor_surreal}}]
  \label{theo:NO-reel}
  The class $\no$ endowed with its ordering $\leq$ and the operations
  $+$, $-$ and $\cdot$ is a totally ordered Field which contains
  $\mathbb{R}$ and \emph{\textbf{On}}.
\end{thm}

$\no$ being an ordered Field carries a \textbf{natural valuation} (see
e.g.\ \cite[Ch. II, Sect. 4, Satz 1]{priess-crampe:book} or
\cite[p.16]{kuhl:ord-exp}) which sends an element to the class of
elements that are \textbf{Archimedean equivalent} -- write $a\asymp b$
-- to itself, i.e.\ for any $a\in\no$, denoting its absolute value by
$|a|=\max\{a,-a\}$:
$$
[a]:=\{b\in\no\, | \, \exists n\in\mathbb{N},\ |a|\leq n \cdot
|b|\textrm{ and } |b|\geq n \cdot |a| \}.
$$
Conway says that $a$ and $b$ are \emph{commensurate}, whereas Gonshor
says that $a$ and $b$ have \emph{same order of magnitude}. Note that
the positive part of such equivalence class is always convex, so it
needs to have a unique element of minimal length. In other words,
there is a \emph{canonical } complete system of representatives of the
Archimedean equivalence classes which forms a cross section of the
value group. Conway finds a way to express naturally such elements as
the images of his so-called \textbf{$\omega$-map}, a map that
generalizes the classical ordinal exponentiation:
\begin{thm}[{\cite[Ch.3, Theorems 19 and 20]{conway_numb-games},
    \cite[Theorems 5.1 to 5.4]{gonshor_surreal}}]
  \label{theo:omega-map}
  The recursive formula:
  $$
  \forall a\in\no, \quad \omega^{a}:=\left\{ 0,\ n \cdot \omega^{a^L}\
    |\ \omega^{a^R}/2^n\right\}
  $$
  (where it is understood that $n$ varies in $\mathbb{N}$) defines an
  ordered Group morphism:
  \begin{align*}
    \Omega: \left(\emph{\no}, +, <\right) & \rightarrow \left(\emph{\no}_{>0}, \cdot, <\right)\\
    a & \mapsto \Omega(a) := \omega^a
  \end{align*}
  that extends the exponentiation with base $\omega$ of the
  ordinals. Moreover, for any $a\in\emph{\no}$, the surreal number
  $\omega^a$ is the positive representative of minimal length in its
  own Archimedean equivalence class.
\end{thm}

Denote by $\omega^{\mathrm{\no}}$ the Image of the $\omega$-map. Note
that the residue field $\mathbb{R}$ of the natural valuation is
naturally embedded in \no. Thus, one has an \emph{approximation
  algorithm}: for any surreal number $a$, there are a surreal number
$b$ and a real number $r$ such that $|a-\omega^b \cdot r|$ is less
than $|a|$ and $|\omega^b \cdot r|$ and in a different Archimedean
class, which we denote as $a-\omega^b \cdot r\prec a$ and
$a-\omega^b \cdot r\prec \omega^b$, and also as
$a \sim \omega^b \cdot r$. Such $b$ is called $\mathrm{Ind}(a)$ in
\cite[Ch.10, Sect.C]{gonshor_surreal}, defining thus a map which is an
incarnation of the natural valuation in the sense of Krull. The value
group is (isomorphic to) $\no$ viewed as an ordered additive Group,
whose multiplicative copy is the cross section
$\omega^{\mathrm{\no}}$. We are in the situation where an adapted
version of \emph{Kaplansky's Embedding Theorem} (see
\cite{kap}\cite[Satz 21, p.62]{priess-crampe:book}) applies. We can
say that $\no$ is \textbf{weakly spherically complete}, namely every
\emph{set} of balls directed by inclusion has non-empty
intersection. Equivalently, any pseudo-Cauchy sequence \emph{indexed
  by some ordinal} has a pseudo limit in $\no$. Subsequently, $\no$ is
isomorphic to the following \textbf{Field of formal power series}
$\mathbb{R}((\omega^{\emph{\no}}))$, i.e.\ the class of formal
expressions $\sum_{i<\lambda}\omega^{a_i} \cdot r_i$ where
$\lambda\in\mathrm{\textbf{On}}$, $(a_i)_{i<\lambda}$ is a strictly
decreasing sequence in \no\ and for any $i<\lambda$,
$r_i\in\mathbb{R}\setminus\{0\}$.

\begin{rem}
  Note that we are abusing the notation
  $\mathbb{R}((\omega^{\emph{\no}}))$, since it usually denotes the
  full \textbf{Field of generalized power series} (also called Hahn
  series field or Malcev-Neumann series field: see
  e.g.\ \cite{matu:gener-series-diff}) which is bigger than
  $\no$. Indeed, in this context, the field of generalized series, as
  a maximally valued field, should consist of series with supports
  being classes. This corresponds to pseudo-Cauchy sequences indexed
  by classes.
\end{rem}

Conway and Gonshor respectively establish two direct
\emph{constructive} proofs of the description of surreal numbers as
power series. They are based on the following key fact: \emph{one does
  have in \emph{\no} a notion of convergence for such generalized
  series}. Indeed, any generalized series
$\sum_{i<\lambda}\omega^{a_i} \cdot r_i$ is defined in $\no$ as the
simplest surreal number having exactly such expansion. In other words,
if $\lambda=\alpha+1$, then
$\sum_{i<\lambda}\omega^{a_i} \cdot r_i:=\sum_{i<\alpha}\omega^{a_i}
\cdot r_i\, +\, \omega^\lambda \cdot r_\lambda$. If $\lambda$ is a
limit ordinal, then $\sum_{i<\lambda}\omega^{a_i} \cdot r_i$ is
defined as the cut obtained as follows for
$\epsilon\in\mathbb{R}_{>0}$:
\begin{multline*}
  \sum_{i<\lambda}\omega^{a_i} \cdot r_i:=\\
  \left\{\sum_{i\leq\alpha}\omega^{a_i} \cdot r_i+\omega^{a_\alpha}
    \cdot (r_\alpha-\epsilon),\, \alpha<\lambda\ |\
    \sum_{i\leq\alpha}\omega^{a_i} \cdot r_i+\omega^{a_\alpha} \cdot
    (r_\alpha+\epsilon),\, \alpha<\lambda\right\}.
\end{multline*}

Consequently, \emph{the approximation algorithm leads to a unique
  expansion of a surreal number as a generalized series.} Indeed, if
the expansion of a surreal number $a$ were indexed over all the
ordinals, then its length $l(a)$ would be greater than \textbf{On}.
Thus we obtain an actual new expression for surreal numbers:

\begin{thm}[Conway normal form of surreal numbers {\cite[Theorems 5.5
    to 5.8]{gonshor_surreal}}]
  \label{theo:NO-series-gene}
  Any surreal number $a\in\emph{\no}$ can be written uniquely as
  $$
  a=\sum_{i<\lambda}\omega^{a_i} \cdot r_i,\ \ \
  \textrm{ the \textbf{normal form} of }a,
  $$
  where $\lambda\in$ \emph{\textbf{On}}, the transfinite sequence
  $(a_i)_{i<\lambda}$ is strictly decreasing and for any $i<\lambda$,
  $r_i\in\mathbb{R}\setminus\{0\}$. Thus:
  $$
  \emph{\no}=\mathbb{R}((\omega^{\no})).
  $$
  where $\omega^{\no}$ is seen as the Group of (generalized)
  monomials.
\end{thm}

We recall that the addition of such generalized series is termwise and
the multiplication is (the straightforward generalization of) the
convolution product for power series (Cauchy product). For ordinal
numbers, the Conway normal form coincides with the classical Cantor
normal form.

\begin{ex}
  \label{ex:analytic}
  For any non-negative integers $n$ and $m$, one has that:
  $$
  \omega^{-n}=\left\{0\, |\, \frac{ \omega^{-(n-1)}}{2^k}\right\},
  $$
  and in particular:
  $$
  \omega^0:=\left\{0\, |\, \emptyset\right\}=1,\qquad\qquad
  \omega^{-1}=\left\{0\, |\, \frac{1}{2^k}\right\}=\frac{1}{\omega} \
  ( k\in\mathbb{N}).
  $$
  Moreover, by induction on $n+m$ and applying
  Theorem~\ref{thm:cofin}:
  \begin{align*}
    \omega^{-n} \cdot \omega^{-m}&= \left\{0\, |\, \frac{ \omega^{-(n-1)}}{2^k}\right\} \cdot \left\{0\, |\, \frac{ \omega^{-(m-1)}}{2^l}\right\}\\
                                 &= \left\{0,\,  \frac{ \omega^{-(m-1)}}{2^l} \cdot{} \omega^{-m}+ \frac{ \omega^{-(n-1)}}{2^k} \cdot \omega^{-m}-  \frac{ \omega^{-(n-1)}}{2^k}  \cdot \frac{ \omega^{-(m-1)}}{2^l}\, | \right. \\
                                 & \qquad \left. \frac{ \omega^{-(n-1)}}{2^k} \cdot \omega^{-m},\, \omega^{-n} \cdot \frac{ \omega^{-(m-1)}}{2^l}\right\}\\
                                 &= \left\{0\, |\, \frac{ \omega^{-(m+n-1)}}{2^i}\right\}\\
                                 &= \omega^{-(n+m)}.
  \end{align*}
  Subsequently, we can compute:
  \begin{align*}
    \sum_{n\geq 0}\omega^{-n}&= \left\{\sum_{n= 0}^N\omega^{-n}\,-\,\omega^{-N} \, |\, \sum_{n= 0}^N\omega^{-n}\,+\,\omega^{-N} \right\}\ (N\in\mathbb{N})\\
                             &=\left\{\frac{1- \omega^{-N-1}}{1- \omega^{-1}}\,-\,\omega^{-N} \, |\,\frac{1- \omega^{-N-1}}{1- \omega^{-1}}\,+\,\omega^{-N}\right\}\\
                             &= \frac{1}{1- \omega^{-1}} \cdot \left\{1-\omega^{-N}\, |\, 1+\omega^{-N}-2 \cdot \omega^{-N-1} \right\}\\
                             &= \frac{1}{1- \omega^{-1}}\,  \cdot \, 1\\
                             &= \frac{1}{1- \omega^{-1}}.
  \end{align*}
\end{ex}

\begin{rem}
  The notion of generalized series goes back to the seminal paper
  \cite{hahn:nichtarchim} of H. Hahn. There, he uses such construction
  of formal power series with real coefficients and exponents in an
  arbitrary ordered \emph{set}, to prove what is called now the
  \emph{Hahn embedding theorem} for ordered Abelian groups. The key
  tool to obtain such an embedding is the Archimedean equivalence
  relation. Hence, the Conway normal form of a surreal number can be
  viewed also as an extension of Hahn's result to \no\ viewed as an
  ordered Abelian Group. Now, considering its multiplicative version
  $\omega^{\mathrm{\no}}$, we can introduce the following natural
  formalism:
  $$
  \omega^a=\omega^{\sum_{i<\lambda}\omega^{a_i} \cdot r_i} =:
  \prod_{i<\lambda}\left(\omega^{\omega^{a_i}}\right)^{r_i}.
  $$
  Note that the numbers $\omega^{\omega^b}$ for $b\in\mathrm{\no}$
  consists of a canonical complete system of representatives of the
  equivalence classes of surreal numbers that are
  \emph{multiplicatively equivalent} in the following sense: for
  $a\in\textrm{\no}_{>0}$, $|a|\succ 1$,
  $$
  [a]_{\textrm{mult}}=[a^{-1}]_{\textrm{mult}} := \{ b^{\pm 1}\in\no
  \, | \, \exists n\in\mathbb{N},\ |a|^n\geq |b|\textrm{ and }
  |b|^n\geq |a| \}.
  $$
\end{rem}

Note that the Conway normal form of any surreal number $a$ can be
split into three parts (each of which possibly trivial) in the
following two manners. Firstly:
\begin{align}
  \label{equ:addit-decompo}
  \begin{split}
    a &= \sum_{a_i>0}\omega^{a_i} \cdot r_i\,+\,r_{i_0} \,+\,\sum_{a_i<0}\omega^{a_i} \cdot r_i\\
    &= \textbf{ purely infinite part } + \textbf{ constant term } +
    \textbf{ infinitesimal part}.
   \end{split}
\end{align}
Secondly:
\begin{align}
  \label{equ:mult-decompo}
  \begin{split}
    a &= \omega^{a_1}\, \cdot \,r_1\, \cdot \,\left(1+\epsilon\right)\\
    & =\textbf{ leading monomial } \cdot \textbf{ leading coefficient
    } \cdot \textbf{ unit}.
  \end{split}
\end{align}
For the sake of notation, let $\mathbb{J} \subseteq \no$ be the
(non-unital) ring of purely infinite surreal numbers.

It follows from the representation of a surreal number as a
generalized power series (with real closed residue field $\mathbb{R}$
and divisible valued group \no) and from Kaplansky's Embedding Theorem
that the Field \no\ is a universal domain for real closed fields in
the following strong sense. A subset $S$ of $\no$ is said to be
\textbf{initial} if it verifies:
$$
\forall (a,b)\in S\times\no,\quad b<_s a\,\Rightarrow\, b\in S.
$$

\begin{thm}[Universal real closed
  Field]
  \label{theo:NO-reel-clos-universel}
  $ $
  \begin{enumerate}
  \item \cite[Ch. 5]{conway_numb-games} \cite[Ch. 5,
    Sect. D]{gonshor_surreal} The proper class \emph{\no} is a real
    closed Field.
  \item \cite[Theorems 28 and 29]{conway_numb-games} \cite[Theorems 9
    and 19]{ehrlich:simpl-hierarc} Any divisible ordered Abelian
    group, respectively any real closed field, is isomorphic to an
    initial subgroup of $(\no,+)$, respectively an initial subfield of
    $(\no,+, \cdot )$.
  \end{enumerate}
\end{thm}

E.g., a surreal number
$a=\omega^{a_1}\, \cdot \,r_1\, \cdot \,\left(1+\epsilon\right)$ has
an $n$th root if and only if its leading coefficient
$r_1\in\mathbb{R}$ has an $n$th root, via the following formula:
$$
a^{1/n}=\omega^{a_1/n}\, \cdot \,r_1^{1/n}\, \cdot
\,\left(1+\epsilon\right)^{1/n}
$$
where $\left(1+\epsilon\right)^{1/n}$ is construed as the
corresponding series expansion in powers of $\epsilon$
(straightforward generalization of the Maclaurin-Taylor series
expansion for the function $\left(1+x\right)^{1/n}$ of a real variable
$x$).

\begin{rem}
  The previous comment as well as Example \ref{ex:analytic} provide a
  glimpse of the following fact \cite[Remark p.43]{conway_numb-games}:
  \emph{any function analytic in some non-empty domain of the real
    plane can be extended to the surreal numbers enclosed in such
    domain via its local power series expansion}. This will also be
  illustrated in the following section. In fact, there is an actual
  theory of analytic functions over surreal numbers
  \cite{alling:analysis-surreal}.
\end{rem}

\subsection{Exponential and logarithmic functions.}

Building on unpublished ideas of Kruskal, Gonshor defines inductively in \cite{gonshor_surreal} a surjective \textbf{exponential map}:
$$\exp:(\no,+,<)\rightarrow(\no_{>0},\, \cdot \,,<)$$ by means of the power series expansion of the real exponential map:
\begin{equation}
  \label{equ:taylor-exp}
  e^x = \sum_{k\geq 0}\frac{x^k}{k!}.
\end{equation}
He obtains also a partial inductive definition of its reciprocal
function, the \textbf{logarithmic map}
$\log:(\no_{>0},\, \cdot \,)\rightarrow(\no,+)$. This terminology is
also justified by the following key fact: these maps coincide with the
usual exponential and logarithm on real numbers, and with (the
extension to surreal numbers of) the Maclaurin-Taylor series
(\ref{equ:taylor-exp}) in powers of $a$ for any infinitesimal surreal
$a$. Let us denote $E_n(x):=\sum_{k= 0}^n\frac{x^k}{k!}$.

\begin{thm}[{\cite[Theorems 10.1 to 10.9, Corollaries 10.1 to
    10.3]{gonshor_surreal}}]
  \label{theo:exp1}
  Consider the recursive formula:
  \begin{gather*}
    \forall a\in\no,\ \exp(a):=\\
    \left\{ 0,\ \exp(a^L)E_n(a-a^L),\ \exp(a^R)E_{2n+1}(a-a^R)\ |\
      \frac{\exp(a^R)}{E_n(a^R-a)},\
      \frac{\exp(a^L)}{E_{2n+1}(a^L-a)}\right\}
  \end{gather*}
  where $n\in\mathbb{N}$ and \emph{in the right hand side, only the
    $n$'s such that $E_{2n+1}(a^L-a)>0$ are considered}. This defines
  an ordered Group isomorphism:
  \begin{align*}
    \exp: \left(\no, +, <\right) & \rightarrow \left(\no_{>0}, \cdot, <\right)\\
    a & \mapsto \exp(a).
  \end{align*}
  The reciprocal isomorphism is denoted:
  \begin{align*}
    \log: \left(\no_{>0}, \cdot, <\right) & \rightarrow \left(\no, +, <\right)\\
    b & \mapsto \log(b).
  \end{align*}
  which verifies for any $b\in\no$:
  \begin{multline*}
    \log(\omega^b):=\left\{ \log(\omega^{b^L})+n,\  \log(\omega^{b^R})-\omega^{(b^R-b)/n}\ |\right.\\
    \left. \log(\omega^{b^R})-n,\
      \log(\omega^{b^L})+\omega^{(b-b^L)/n}\right\}.
  \end{multline*}

  Moreover, for any $r\in\mathbb{R}$, $s\in\mathbb{R}_{>0}$ and
  $a\in\no$, $a\prec 1$:
  $$
  \left\{
    \begin{aligned}
      \exp(r) & = e^r\\
      \exp(a) & = \sum_{k\geq 0}\frac{a^k}{k!}
    \end{aligned}
  \right.
  \quad \textrm{ and } \quad
  \left\{
    \begin{aligned}
      \log(s) & = \ln(s)\\
      \log(1+a) & = \sum_{k\geq 1}\frac{(-1)^{k+1}a^k}{k},
    \end{aligned}
  \right.
  $$

  $$
  \left\{
    \begin{aligned}
           \exp(r)=e^r\\
           \exp(a)=\sum_{k\geq 0}\frac{a^k}{k!}
    \end{aligned}
  \right.
  \quad \textrm{ and } \quad
  \left\{
    \begin{aligned}
      \log(s) & = \ln(s)\\
      \log(1+a) & = \sum_{k\geq 1}\frac{(-1)^{k+1}a^k}{k}.
    \end{aligned}
  \right.
  $$
\end{thm}

Let us give a hint on how to verify the validity of the inductive
formula for exp. Since $E_n$ is a truncation of the full series
(\ref{equ:taylor-exp}), one can observe that the elements involved in
the cut verify for instance:
\begin{multline*}
  \exp(a^L)E_n(a-a^L) < \exp(a^L) \cdot \exp(a-a^L) = \exp(a) = \\
  = \frac{\exp(a^L)}{\exp(a^L-a)} < \frac{\exp(a^L)}{E_{2n+1}(a^L-a)}
\end{multline*}
whenever only the $n$'s for which $E_{2n+1}(a^L-a)>0$ are considered.

\medskip

The computation of the exponential for \textbf{finite} surreal numbers
(i.e.\ of type $x=r+\epsilon$ for some real $r$ and some infinitesimal
$\epsilon$) amounts to determining the sum of their power series via
(\ref{equ:taylor-exp}). For purely infinite surreal numbers (i.e.\
numbers whose Conway normal form has only positive exponents), Gonshor
obtains specific results. As an example, let us compute:
\begin{align*}
  \exp(\omega) &= \exp\left(\left\{n\ |\ \emptyset\right\}\right)\\
               &= \left\{ 0,\ \exp(n) \cdot E_n(\omega-n)\ |\  \emptyset\right\} \ (n\in\mathbb{N})\\
               & \qquad\qquad \textrm{ (since }E_{2n+1}(n-\omega)<0\textrm{ for any } n)\\
               &= \left\{\omega^n\ |\ \emptyset\right\}\qquad \textrm{ (by cofinality)}\\
               &= \omega^\omega.
\end{align*}
Similarly, one can show that for some of the subsequent infinite
ordinals, exp coincides with the omega map, which itself coincides
with ordinal exponentiation. However, this is not a general rule. For
instance, let us consider the first \emph{epsilon number}
$\varepsilon_0=\{1,\, \omega,\, \omega^\omega,\ldots\ |\
\emptyset\}=\{\Omega^n(0)\ |\ \emptyset\}$. Recall that this is the
first fixed point of the ordinal exponentiation:
$\omega^{\varepsilon_0}=\varepsilon_0$. Supposing that we have already
proved that $\exp$ coincides with $\Omega$ on the infinite ordinals
$\Omega^n(0)$, $n>1$, we compute as before:
\begin{align*}
  \exp(\varepsilon_0) &= \exp\left(\left\{\Omega^n(0)\ |\ \emptyset\right\}\right) \ \  \ \ \ \  \ (n\in\mathbb{N})\\
                      &= \left\{ 0,\ \exp\left(\Omega^n(0)\right) \cdot E_m\left(\varepsilon_0-\Omega^n(0)\right)\ |\  \emptyset\right\}\ (m\in\mathbb{N})\\
                      &= \left\{\Omega^{n+1}(0) \cdot \left(\varepsilon_0-\Omega^n(0)\right)^m\ |\ \emptyset\right\}\ \ \  \ \  \textrm{ (by cofinality)}\\
                      &= \left\{\omega^{\varepsilon_0 \cdot m}\ |\ \emptyset\right\}\qquad  \textrm{ (since }\omega^{\varepsilon_0}=\varepsilon_0 \textrm{ and again by cofinality)}\\
                      &= \omega^{\varepsilon_0 \cdot \omega} \,=\, \omega^{\omega^{\varepsilon_0+1}}\,\neq\, \varepsilon_0.
\end{align*}
The following theorem describes the exact relation between the
exponential map and the omega map for purely infinite surreal numbers.

\begin{thm}[{\cite[Theorems 10.8 to 10.13]{gonshor_surreal}}]
  \label{theo:exp2}
  Monomials correspond to exponentials of purely infinite surreal
  numbers:
  $$
  \omega^{\no}=\exp( \mathbb{J}).
  $$
  Moreover for any purely infinite surreal number
  $a=\sum_{i<\lambda}\omega^{a_i}r_i$ (i.e.\ $\forall i,\, a_i>0$):
  $$
  \exp(a)= \omega^{y} \textrm{ with }
  y=\sum_{i<\lambda}\omega^{g(a_i)}r_i
  $$
  where the map $g:\no_{>0}\rightarrow \no$ is defined by:
  $$
  g(a):=\left\{ \mathrm{Ind}(a),\ g(a^L)\ |\ g(a^R)\right\}.
  $$
  The inverse map $\log$ of $\exp$ satisfies for any surreal number
  $b=\sum_{j<\lambda}\omega^{b_{j}}s_{j}$:
  $$
  \log(\omega^b)=\log(\prod_{j<\lambda}
  \left(\omega^{\omega^{b_{j}}}\right)^{s_{j}})=
  \sum_{j<\lambda}\omega^{h(b_j)}
  $$
  where the map $h:\no\rightarrow \no_{>0}$ is defined by:
  $$
  h(b):=\left\{ 0,\ h(b^L)\ |\ h(b^R),\ \omega^b/2^n\right\}.
  $$
  We have $h=g^{-1}$ on $\no$.
\end{thm}

Subsequently, Gonshor gives detailed results on the map $g$. For
instance, he proves in \cite[Theorem 10.17]{gonshor_surreal} that
$g(a)=a$ -- equivalently, exp coincides with $\Omega$ on $\omega^a$ --
for any surreal number $a$ such that
$\varepsilon_\alpha+\omega\leq a\leq \beta<\varepsilon_{\alpha+1}$ or
such that $1\leq a\leq \beta< \varepsilon_0$ for some
$\alpha,\beta\in\textbf{\textrm{On}}$ ($\varepsilon_\alpha$ denotes as
usual the $\alpha$'th epsilon number). In turn, Gonshor shows that $g$
fails to be the identity map in a close neighborhood of any of the
\emph{generalized epsilon numbers}, i.e.\ the fixed points for the
$\omega$-map \cite[Ch.9]{gonshor_surreal}. This illustrates a key
difference between $\Omega$ and exp, since the latter cannot have
fixed points.

\medskip

As it has been underlined in \cite{ehrlich-vdd:surreal-exp}, the
results on $\no$ described so far agree with important results on the
model theory of the ordered field of real numbers with restricted
analytic functions and the exponential function. More precisely, let
$L_{\mathrm{an}}(\exp)$ be the language of ordered rings augmented by
a symbol for each real multivariate power series convergent on the
closed unit hypercube corresponding to their domain, and by a symbol
exp, and let $T_{\mathrm{an},\exp}$ be the
$L_{\mathrm{an}}(\exp)$-theory of the field of real numbers where the
$L_{\mathrm{an}}(\exp)$-symbols are naturally construed as:
\begin{itemize}
\item the real analytic functions corresponding to their convergent
  power series expansion on the hypercube on and extended by 0 outside
  of it to a well-defined map on $\mathbb{R}$;
\item the real exponential map for exp, with the axiomatization proved
  by Ressayre \cite{ress}.
\end{itemize}
The corresponding structure $\mathbb{R}_{\mathrm{an},\exp}$ is
\emph{model complete}, \emph{o-minimal} and has a \emph{complete
  axiomatization} with $T_{\mathrm{an},\exp}$ \cite[Second Main
Theorem]{wil:modelcomp} \cite[Corollaries 4.5, 4.6 and
5.13]{dmm:exp}. Recall that these key model theoretic properties have
striking geometric counterparts, which are the so-called
\emph{tameness properties} after Grothendieck's
\cite{grothendieck:esquisse}: finiteness of the number of connected
components, cell decomposition, dimension theory,
triangulation,... See \cite{vdd:tame} and \cite{coste:intro_ominimal}.
Using the previously described interpretation of
$L_{\mathrm{an}}(\exp)$ and the corresponding results in Gonshor's
work, the authors obtain:

\begin{thm}[{\cite[Theorem 2.1]{ehrlich-vdd:surreal-exp}}]
  The structure $\emph{\no}_{\mathrm{an},\exp}$ is an elementary
  extension of $\mathbb{R}_{\mathrm{an},\exp}$.
\end{thm}

Moreover, the authors describe a family of initial subfields of \no,
each of which providing an elementary substructure of
$\emph{\no}_{\mathrm{an},\exp}$ as well as an elementary extension of
$\mathbb{R}_{\mathrm{an},\exp}$. This means that all these fields
share the same first order properties, in particular \emph{the
  geometry based on these non-Archimedean fields is tame}.

\section{From Hardy fields to transseries and $H$-fields}
\label{sec:H-fields}

\subsection{Hardy fields and tame geometry.}
\label{sec:hardy-fields}

Few years after Gonshor's \cite{gonshor_surreal}, several other
similar non standard models of $\mathbb{R}_{\mathrm{an},\exp}$
appeared in two different areas in connection with new important
results:
\begin{itemize}
\item in model theory, the field of Logarithmic-Exponential series
  \cite{vdd:LE-pow-series,dmm:exp} and the fields of
  Exponential-Logarithmic series \cite{kuhl:ord-exp,ks:kappa-bdd-exp}
  related to the results already cited on the theories of
  $\mathbb{R}_{\exp}$ and $\mathbb{R}_{\mathrm{an},\exp}$;
\item in differential equations, fields of grid-based and well-ordered
  transseries \cite{vdh:autom-asymp,vdh:transs_diff_alg,schm01}
  related to the proof of the Dulac conjecture \cite{ecalle:dulac}.
\end{itemize}

In any case, one of the main purposes of these objects is to provide
formal analogues to non oscillating functions, in the sense of
\emph{asymptotic scales}. In terms of real functions, exp and log have
been used at an earlier stage to describe asymptotically non
oscillating functions. Many interesting functions have an asymptotic
behavior at $+\infty$ that can be described in terms of $\exp$, $\log$
and algebraic functions -- for instance, the prime number theorem says
that $\pi(x) \sim \frac{x}{\log(x)}$, or Stirling's approximation
states that $n! = \Gamma(n) \sim \sqrt{2\pi
  n}e^{n\log(n)-n}$. Especially, one can obtain asymptotic equivalents
at $+\infty$ for solutions of differential equations (Airy functions,
hypergeometric functions, etc.).

The theory of such asymptotic expression was eventually formalized by
Hardy in his \emph{logarithmico-exponential functions} (LE-functions)
\cite{hardy}, after the work of Du Bois-Reymond, with further
inspiration from Liouville's work on elementary functions. We
paraphrase Hardy's definition using germs instead of actual
functions. We define the \textbf{germ} (at $+\infty$) of a real
function $f : (a, +\infty) \to \mathbb{R}$ as the equivalence class:
$$
  [f] := \{ g : (b, +\infty) \to \mathbb{R} \,:\,\exists c \in
  \mathbb{R}\ \forall x > c\ (f(x) = g(x)) \}.
$$
Germs of real functions form a ring -- let us denote it by
$\mathcal{G}$ -- with the identity given by the constant function
$[1]$ and the operations $[f] + [g] = [f + g]$ and
$[f \cdot g] = [f] \cdot [g]$.

By definition, the class of \textbf{logarithmico-exponential germs}
(LE-germs) is the smallest ring of germs such that:
\begin{itemize}
\item the germs of the constant functions and of the identity function
  $x : \mathbb{R} \to \mathbb{R}$ are in the ring;
\item if $[f]$ is in the class, then $[\exp(f)]$ is in the ring;
\item if $[f]$ is in the class and the image of $f$ is strictly
  positive, then $[\log(f)]$ is in the ring;
\item if $[f]$ is in the class, and the image of $f$ is contained in
  the domain of a semi-algebraic function $A$, $[A(f)]$ is in the
  ring.
\end{itemize}
Note that, by construction, any LE-function $f$ is differentiable in a
neighborhood of $+\infty$ and its derivative $f'$ is again a
LE-function. The mapping $[f]\mapsto[f']$ is a well-defined derivation
on LE-germs. Moreover, as a crucial result \cite[Theorem,
p.18]{hardy}, LE-functions are eventually continuous, of constant sign
and monotonic. Subsequently, LE-germs form an ordered field, and
eventually a \emph{differential ordered field}. Therefore, one can
compare the order of growth of logarithmico-exponential germs. We
define (compare with the definitions given before and after Theorem
\ref{theo:omega-map}):
\begin{itemize}
\item $[f] \prec [g]$ if
  $\lim_{x \to +\infty} \displaystyle\frac{f}{g}(x) = 0$;
  equivalently, if $n \cdot [f] < [g]$ for all $n \in \mathbb{N}$;
\item $[f] \preccurlyeq [g]$ if $[g] \not\prec [f]$;
\item $[f] \asymp [g]$ if $[f] \preccurlyeq [g]$ and
  $[g] \preccurlyeq [f]$;
\item $[f] \sim [g]$ if $[f] - [g] \prec [f]$ and
  $[f] - [g] \prec [g]$.
\end{itemize}
For instance, we have $[e^x] \succ [x^n]$ for all $n \in \mathbb{N}$,
while $[e^x] \sim [e^x] + [\log(x)]$, and so on. In fact, for any two
germs $[f]$, $[g]$, one has either $[f] \prec [g]$, $[g] \prec [f]$,
or $[f] \sim [rg]$ for some non-zero real number $r \in \mathbb{R}$.

LE-germs represent a large asymptotic scale, which includes the usual
scale consisting of powers of $x$, and that has still the rich
algebraic structure of an ordered differential field. Building on
Hardy's work, Bourbaki introduced the following general notion:

\begin{defn}[{\cite[p.V-36]{bour:unevar}}]
  \label{defi:hardy-field}
  A \textbf{Hardy field} is a subring of $\mathcal{G}$ that is a field
  and that is closed under differentiation.
\end{defn}

For instance, the field of (germs of) rational fractions with real
coefficients $\mathbb{R}(x)$ is a Hardy field. So is the field of
functions real meromorphic at infinity $\mathbb{R}\{1/x\}$. There are
various ways to derive new Hardy fields from a given one $K$: by
taking its real closure $K^{rc}$ in $\mathcal{G}$
\cite{robinson:rc-hardy-field}; by taking its \textbf{Liouville
  closure} (i.e.\ closing under integration and logarithmic
integration) $K^{lc}$ in $\mathcal{G}$ \cite{bour:unevar}; by taking
its closure under resolution of differential equations $Q(y)y'=P(y)$
for $P,Q\in K[y]$ \cite[based on an unpublished work by M.\
Singer]{rosenlicht:hardy_fields}; by adjoining solutions of certain
second order linear differential equations
\cite{boshernitzan:2nd-order-ode,ros:order2-ODE}. Note that the field
of LE-germs is not closed under the latter extensions; by Liouville's
theorem, it does not contain the antiderivative of $e^{-x^2}$.

In fact, Hardy fields tends to be central objects in tame geometry:
\emph{given an o-minimal structure over $\mathbb{R}$, the germs at
  $+\infty$ of unary real functions that are definable (i.e.\ whose
  graph is a subset definable in the structure) form a Hardy
  field}. Indeed, differentiation is a definable operation, and
functions are ultimately non oscillating by o-minimality. The simplest
examples are the Hardy fields corresponding to the structure of
\emph{semialgebraic sets} and to \emph{subanalytic sets}: they consist
in Puiseux power series at $+\infty$ with real coefficients and that
are algebraic, respectively convergent. Another key example, which
contains both the preceding examples and the field of LE-germs, is the
Hardy field corresponding to $\mathbb{R}_{\mathrm{an},\exp}$. Note
that, in this case also, one has an explicit description of these
germs in terms of compositions of exp, log and restricted real
analytic functions \cite[Corollary 4.7]{dmm:exp}.

Nevertheless, understanding Hardy fields in full generality remains a
challenging problem. Generally speaking, the union of two Hardy fields
does not generate a Hardy field. It implies that there is no maximum
Hardy field, but several maximal ones. Furthermore, there exist Hardy
fields that contain functions growing more than any iteration of the
exponential function \cite{bosher:hardy}. Several questions remain
unanswered: how do maximal Hardy fields look like? What differential
equations can be solved in them? Are there maximal Hardy fields
corresponding to some o-minimal structure?

\subsection{From the standpoint of formal power series.}
\label{sec:transseries}

As we already announced, there are several similar non standard models
of the theory of $\mathbb{R}_{\mathrm{an},\exp}$ based on generalized
power series: fields of LE-series, EL-series and transseries. All of
them are special kinds of ordered fields of formal power series, which
carry exponential and logarithmic maps. Recall that, according to
\cite{kks:exp-power-series}, there is no hope for defining a global
exponential function directly on a full field of generalized power
series $\R((\mathfrak{M}))$ (where $\mathfrak{M}$ denotes the
(multiplicatively written) ordered Abelian group of (generalized)
monomials). Indeed, such exp should give an isomorphism between the
additive group and the positive part of the multiplicative group. By
decomposing the additive group as in \eqref{equ:addit-decompo} and the
multiplicative group as in \eqref{equ:mult-decompo}, one can verify
that the isomorphism induces three isomorphisms between the three
components of the decompositions. However, even though such
isomorphisms exist in any generalized series field on the elements of
non negative valuation (i.e., the ``finite series'', compare with
Theorem \ref{theo:exp1}), there cannot be any between the additive
group of series having terms with negative value (i.e.\ ``purely
infinite'' series) and the multiplicative group of monomials
$\mathfrak{M}$, as proven in \cite{kks:exp-power-series}. Thus, no
global exponential function exist.

Hence, one has to use a special construction to obtain such
isomorphism in the context of generalized power series. The main
common idea in the different versions is due independently to
\cite{dahn:limit-exp-terms,dahn-gorhing} and \cite{ecalle:dulac}, and
may be seen as an abstract version of Hardy's construction of log-exp
functions \cite{hardy}. It consists in building extensions of a given
$\R((\mathfrak{M}))$ which will contain the missing elements to obtain
a well-defined exp and a surjective log: the corresponding field with
global exp and log is then obtained as a subfield of such an enlarged
field of generalized series. Let us describe the process for the
\textbf{exponential closure} of $\R((\mathfrak{M}))$:
\begin{itemize}
\item Any field of generalized series $\R((\mathfrak{M}))$ has a
  partial exponential map $e$ defined on the series with non negative
  value (the ``finite part'' of the series) using the real exponential
  map and the Taylor formula (\ref{equ:taylor-exp}).
\item Given a field of generalized series $(\R((\mathfrak{M})),\, e)$
  endowed with a partial exponential map $e$, we consider
  $\mathfrak{M}^\sharp =: e^\sharp
  \left(\R((\mathfrak{M}_{\succ1}))\right)$ a multiplicative copy of
  the group of purely infinite series where the morphism $e^\sharp$ is
  intended to extend $e$. More precisely, the image of $e$ and its
  complement are ordered lexicographically. Accordingly, the
  \textbf{exponential extension} of $(\R((\mathfrak{M})),\, e)$ is the
  field $(\R((\mathfrak{M}^\sharp)),e^\sharp)$ for $e^\sharp$ the
  corresponding function extension of $e$.
\item Considering infinitely many iterated exponential extensions of
  the given series field $\R((\mathfrak{M}))$, its exponential closure
  -- denote it by $\R((\mathfrak{M}))^{\mathrm{E}}$ -- is the direct
  limit of such tower of fields with partial exponentials. The
  morphism obtained as a limit of the partial $e$'s is a surjective
  exponential that we denote by $\exp$.
\end{itemize}
Note that $\R((\mathfrak{M}))^{\mathrm{E}}$ is a subfield of the
generalized series field $\R((\mathfrak{M}^{\mathrm{E}}))$ where
$\mathfrak{M}^{\mathrm{E}}$ denotes the union of the iterated extended
groups of monomials. Concerning the logarithmic closure, different
approaches provide different objects: see \cite{kuhl-tressl:el-le} for
a comparison between LE-series and EL-series. There are also variants
in case one imposes restrictions on the supports of the powers series
considered: compare grid-based and well-ordered transseries in
\cite{vdh:autom-asymp}. Note also that restricted analytic functions
extend to such various fields as they extend to generalized power
series fields (direct consequence of Neumann's Lemma
\cite{neumann:ord-div-rings}).

\medskip

A key additional feature of these non standard models of
$\mathbb{R}_{\mathrm{an},\exp}$ is that they can be endowed with
particular nice derivations (provided that the initial field
$\R((\mathfrak{M}))$ carries such a derivation):
\begin{itemize}
\item that agree with the structure of power series, in particular
  which commute with infinite sums \cite[Section 3]{dmm:LE-series}
  \cite[Definition 4.1.1]{schm01} \cite[Definition
  3.6]{matu-kuhlm:hardy-deriv-EL-series};
\item that are compatible with the exp and log functions, namely
  $D(\exp(x))=\exp(x) \cdot D(x)$ and $x \cdot D(\log(x))=D(x)$;
\item that behave very like the derivation of a Hardy field, in
  particular like the Hardy field associated to
  $\mathbb{R}_{\mathrm{an},\exp}$ (see \cite[Section 4]{dmm:LE-series}
  \cite[Definition 4.1]{matu-kuhlm:hardy-deriv-EL-series} and the
  following section).
\end{itemize}

Recall that these power series with exp and log are meant to be formal
counterparts to actual non oscillating functions, in particular in the
context of differential equations. This will be discussed further
later on. Let us just state already the following striking embedding
theorem which follows directly from the cited model theoretic results
on $\mathbb{R}_{\mathrm{an},\exp}$:

\begin{thm}[{\cite[Corollary 3.12]{dmm:LE-series}}]
  The Hardy field associated to $\mathbb{R}_{\mathrm{an},\exp}$ embeds
  naturally as a differential exponential and logarithmic real closed
  field into the field of LE-series.
\end{thm}

In \cite{schm01}, the author proposes an axiomatic version of
transseries. In order to avoid confusion, we slightly modify the
terminology used by this author:

\begin{defn}
  \label{defi:transseries}
  Let $\R((\mathfrak{M}))$ be a field of generalized series and $\log$
  be a function such that:
  \begin{description}
  \item[(T1)] the domain of $\log$ consists of the positive series;
  \item[(T2)]
    $\log(\mathfrak{M})\subseteq \R((\mathfrak{M}_{\succ 1}))$;
  \item[(T3)]
    $\log(1+\varepsilon)=\displaystyle\sum_{n\geq
      1}(-1)^{n+1}\displaystyle\frac{\varepsilon^n}{n}$ for any
    $\varepsilon\in\R((\mathfrak{M}_{\preccurlyeq 1}))$;
  \item[(T4)] let $(\mathfrak{m}_n)_{n\in\N}$ be a sequence of
    monomials with $\mathfrak{m}_{n+1}$ in the support of
    $\log(\mathfrak{m}_n)$ for all $n$. There is an integer $n_0$ such
    that for any $n\geq n_0$, for any $\mathfrak{m}$ in the support of
    $\log(\mathfrak{m}_n)$, one has that
    $\mathfrak{m}\succcurlyeq \mathfrak{m}_{n+1}$ and the coefficient
    of $ \mathfrak{m}_{n+1}$ in $\log(\mathfrak{m}_n)$ is $\pm1$.
  \end{description}
\end{defn}

Note that such axioms remain verified by the corresponding transseries
field $\R((\mathfrak{M}))^{\mathrm{E}}$. In fact, the author does even
consider \emph{transfinite} exponential extensions. Axioms (T1) to
(T3) mean that $\R((\mathfrak{M}))$ is endowed with a so-called
\emph{prelogarithm} \cite{kuhl:ord-exp}. (T4) takes into account the
possible existence of elements of type:
$$
y = \exp(x + \exp(\log_2(x) + \exp(\log_4(x) + \dots))),
$$
the latter being viewed as a formal solution to the functional
equation: $y = \exp(x + \log_{2}(y))$ ($x$ is construed as the germ of
identity at $+\infty$). Nevertheless, (T4) forbids the existence of
series of type:
$$
\exp(x + \exp(\log_2(x) + \exp(\log_4(x) + \dots + \log_6(x)) +
\log_4(x)) + \log_2(x)).
$$
The reader should compare this with the notion of \emph{irreducible
  surreal numbers} in \cite[p. 33]{conway_numb-games}. In
\cite{matu-kuhlm:surreel-transseries-EL}, the authors introduced the
more restrictive notion of \textbf{ field of exp-log transseries} as a
common axiomatic for EL-series and transseries, where (T4) is replaced
by:
\begin{description}
\item[(ELT4)] let $(\mathfrak{m}_n)_{n\in\N}$ be a sequence of
  monomials with $\mathfrak{m}_{n+1}$ in the support of
  $\log(\mathfrak{m}_n)$ for all $n$. There is an integer $n_0$ such
  that $\log(\mathfrak{m}_{n_0+k})=\mathfrak{m}_{n_0+k+1}$ for any
  $k\in\mathbb{N}$.
\end{description}
As will be discussed with further details in Section
\ref{sec:derivations} this axiom applies to an important strict
subfield of \textbf{No}, but not to \textbf{No} itself.

\subsection{$H$-fields as a common axiomatic framework}

Concerning Hardy fields as described in Section
\ref{sec:hardy-fields}, Rosenlicht studied the remarkable properties
linking the derivation and the ordering as a \emph{differential
  ordered fields} -- or accordingly the corresponding natural
valuation $v$ as a \emph{differential valued fields}
\cite{rosenlicht:rank}. Building on his work, M. Aschenbrenner and
L. van den Dries introduced in \cite{ad:H-fields} the following more
abstract notion:
\begin{defn}
  An \textbf{$H$-field} is an ordered field $(K, +, \cdot, \leq)$
  equipped with a function $D : K \to K$ such that:
  \begin{enumerate}
  \item additivity: $D(x + y) = D(x) + D(y)$;
  \item Leibniz rule: $D(xy)=xD(y)+yD(x)$;
  \item H1: if $x > \ker(D)$, then $D(x)>0$;
  \item H2: if $|x| < c$ for some $c \in \ker(D)$, then there exists
    $d \in \ker(D)$ such that $|x - d| < c$ for all $c \in \ker(D)$.
  \end{enumerate}
\end{defn}
For a Hardy field $K$, Condition H1 follows from the fact that if
$[f] > \mathbb{R}$, then $\lim_{x \to +\infty}f(x) = +\infty$, which
implies that $f'(x)$ is eventually positive, so $[f'] > 0$. Condition
H2 holds whenever
$\ker \left(\frac{d}{dx}\right)=\mathbb{R}\subseteq K$.

Given an $H$-field $K$ with derivation $D$, the convex hull of
$\ker(D)$ is a valuation ring. Let us denote by
$v : K^{\times} \to \Gamma$ the corresponding valuation, where
$\Gamma$ is a suitable ordered Abelian group ($v(x) \geq 0$ if and
only if $x$ is in the convex hull of $\ker(D)$). We shall use again
the notation of Hardy and Du Bois-Reymond accordingly to this
valuation (compare with Section \ref{sec:hardy-fields}). For
$x, y \in K^{\times}$, we shall write:
\begin{itemize}
\item $x \prec y$ if $v(x) > v(y)$;
\item $x \preceq y$ if $v(x) \geq v(y)$;
\item $x \asymp y$ if $v(x) = v(y)$;
\item $x \sim y$ if $v(x - y) > v(x)$ and $v(x - y) > v(y)$.
\end{itemize}

Derivation and valuation are linked via specific properties. For the
sake of notation, we shall denote by $LD(x)$ the \textbf{logarithmic
  derivative} $LD(x) := D(x)/x$. In particular, one has:
\begin{itemize}
\item a strong version of \textbf{L'Hospital rule}:
  $\forall x,y \not \asymp 1,\ (x \prec y \Leftrightarrow D(x) \prec
  D(y))$;
\item a rule for $LD$:
  $\forall x\prec y\prec 1,\ LD(x)\succcurlyeq LD(y)$.
\end{itemize}
Note that these properties were also used in
\cite{matu-kuhlm:hardy-deriv-EL-series,
  matu-kuhlm:hardy-deriv-gener-series} by the authors as axioms of a
so-called \emph{Hardy type derivation} for a differential valued
field.

In most cases, one also requires that the derivation preserves
infinitesimal elements.

\begin{defn}
  An $H$-field $(K, D)$ has \textbf{small derivation} if
  $D(x) \prec 1$ for all $x \in K$ such that $x \prec 1$.
\end{defn}

Note that this assumption is rather harmless. Indeed, for elements
$y,z \in K$, $z \succ 1$, such that $y = D(z)$, if we set
$D^y : x \mapsto y \cdot D(x)$, then $(K, D^y)$ is an $H$-field with
small derivation. Note also that all Hardy fields with their natural
derivation $\frac{\mathrm{d}}{\mathrm{d}x}$ are $H$-fields with small
derivation ($[f]$ infinitesimal means that
$\lim_{x \to +\infty}f(x) = 0$, which implies that
$\lim_{x \to +\infty}f'(x) = 0$, so $[f']$ is infinitesimal as well).

Rosenlicht observed in \cite{rosenlicht:value-group} that some of the
key relations between the derivation $D$ and the valuation $v$ can be
encoded within the valued group $\Gamma$ equipped with an extra
function that we denote also by $LD$: for
$\gamma = v(a) \in \Gamma^{\neq 0}$, let
$LD(\gamma) := v(LD(a)) = v(D(a)/a) = v(D(a)) - \gamma$. Such
$(\Gamma,LD)$ is called an \textbf{asymptotic couple} in the context
of Hardy fields. In \cite{ad:H-fields}, the authors enhance this study
of asymptotic couples for $H$-fields. In particular, of prime interest
is the following notion related to the question of integration:

\begin{defn}
  Let $(K, D)$ be an $H$-field. Given some $a \in K$, an element
  $b \in K$ is an \textbf{asymptotic integral} of $a$ if
  $a \sim D(b)$.
\end{defn}

\begin{thm}[\cite{rosenlicht:rank}]
  \label{thm:asymp-integr}
  Let $(K, D)$ be an $H$-field, with asymptotic couple $(\Gamma,
  LD)$. For all $a \in K$, $v(a)$ is not the supremum of
  $LD(\Gamma^{\neq 0})$ if and only if $a$ admits an asymptotic
  integral in $K$.
\end{thm}

\begin{cor}
  Let $(K, D)$ be an $H$-field, with asymptotic couple $(\Gamma,
  LD)$. If $LD(\Gamma^{\neq 0})$ has no supremum in $\Gamma$, then
  every element of $K$ admits an asymptotic integral.
\end{cor}

When all elements of $K$ admits an asymptotic integral, we say that
$K$ admits asymptotic integration. This notion applies particularly in
the context of \emph{spherically complete fields}, the abstract
version of fields of generalized power series (see Section
\ref{sec:universal-RCF}).

\begin{thm}[\cite{fvk:ultrametric-spaces}]
  \label{thm:integr-spher-complete}
  Let $K$ be an $H$-field. If every $a \in K$ admits an asymptotic
  integral, and $K$ is spherically complete, then every $a \in K$
  admits an integral.
\end{thm}

Indeed, if one wants to find an integral of a given element $a \in K$,
one can define a sequence of approximations as follows: first, we let
$b_0$ be an asymptotic integral of $a$, then we define inductively
$b_{i+1}$ to be the asymptotic integral of $D(b_i) - a$. Provided one
takes appropriate precautions, such as using H2 to ensure that
$v(b_i) \neq 0$ for all $i$, the resulting sequence $(b_i)$ is
pseudo-Cauchy, and when $K$ is spherically complete, it has a
pseudo-limit $b_{\omega}$. If $D(b_{\omega}) = a$, we are done. If
not, we keep iterating on all ordinals. The procedure eventually
produces an integral for $a$.

\medskip

As mentioned before, Hardy fields can be closed under algebraic
extensions and application of $\exp$ and $\log$ via the notion of
Liouville closure (see Section \ref{sec:hardy-fields}). Similarly,
every $H$-field can be extended to a \textbf{Liouville-closed}
$H$-field, that is to say, a real closed $H$-field in which the
equations $y' = a$ and $z' = az \wedge z \neq 0$ have solutions for
all $a$. Moreover, suppose now that the $H$-field $K$ is also equipped
with an exponential function $\exp : K \to K^{> 0}$ compatible with
the derivation $D$, namely such that $D(\exp(x)) = \exp(x)D(x)$ for
all $x \in K^{> 0}$. Then $K$ is Liouville-closed if and only if the
equation $y' = a$ has a solution for all $a \in K$. Indeed, to solve
$z' = az \wedge z \neq 0$ it suffices to take a $b \in K$ such that
$D(b) = a$ and note that $\exp(b)$ is a solution to $z' = az$.

\begin{thm}[\cite{ad:H-fields}]
  Every $H$-field has one or two Liouville closures up to isomorphism.
\end{thm}

In many important cases, the number of Liouville closures of an
$H$-field $K$ is determined by the structure of its asymptotic couple
(however, A.\ Gehret proved recently that when $K$ has asymptotic
integration, the determining factor for the number of Liouville
closures is not the asymptotic couple, but rather the property of
``$\lambda$-freeness'' of $K$ \cite{gehret}). In
\cite{adh:mod-theo-transseries}, the authors identified a key first
order axiom that applies to the field of LE-series and which implies
uniqueness of the Liouville closure.

\begin{defn}
  An $H$-field $(K, D)$ is \textbf{$\upomega$-free} if it is real
  closed, it admits asymptotic integration, and
  $$
  \forall a\, \exists b\, (b \succ 1 \wedge a - \omega(LD(LD(b)))
  \succeq LD(b)^2.
  $$
  where $\omega(y):=-(2y'-y^2)$.
\end{defn}

$\upomega$-freeness relates to the solvability within $K$ of second
order linear equations of type: $4y''+fy=0$ which is governed by the
behaviour of the corresponding Riccati operator $\omega$. This is
already discussed in \cite{ros:order2-ODE} in the context of Hardy
fields.

\begin{thm}[{\cite[Corollary 13.6.2]{adh:mod-theo-transseries}}]
  If an $H$-field $K$ is $\upomega$-free, then it has a unique
  Liouville closure up to isomorphism.
\end{thm}

We conclude this section by quoting one of the main theorems of
Aschenbrenner, van den Dries and van der Hoeven in
\cite{adh:mod-theo-transseries}, after having been conjectured by the
same authors 20 years ago. The other key notion is the one of
\textbf{Newtonianity}, a version with asymptotic constraints of
differential Henselianity: any asymptotic algebraic differential
equation of Newton degree one has a root in the valuation ring.

\begin{thm}[{\cite{adh:mod-theo-transseries}}]
  \label{thm:model-complete}
  Let $T$ be the first-order theory in the language
  $\mathcal{L} = \{0, 1, +, \cdot, <, \prec, D\}$ whose models are the
  $H$-fields $K$ such that:
  \begin{itemize}
  \item the derivation is small;
  \item $K$ is Liouville-closed;
  \item $K$ is $\upomega$-free;
  \item $K$ is Newtonian.
  \end{itemize}
  Then $T$ is complete and model-complete. Moreover, the fields of
  $LE$-series and grid-based transseries are a model of $T$.
\end{thm}

We can now state the relevant related questions: can $\no$ be given
the structure of an $H$-field? Is it a model of the theory $T$? Is it
universal?

\section{Surreal derivations}
\label{sec:derivations}

Inspired by the definition of $H$-field, and by the properties of
derivations on fields of transseries, the authors define in
\cite{ber-man:surreal-derivations} derivations on $\no$ as follows.

\begin{defn}\label{def:surreal-derivation}
  A \textbf{surreal derivation} is a function $D : \no \to \no$
  satisfying the following properties:
  \begin{enumerate}
  \item Leibniz rule: $D(xy)=xD(y)+yD(x)$;
  \item strong additivity:
    $D \left(\sum_{i\in I} x_{i} \right) = \sum_{i\in I} D(x_{i})$ if
    $(x_{i} \,:\, i \in I)$ is summable;
  \item compatibility with exponentiation: $D(\exp(x)) = \exp(x)D(x)$;
  \item constant field $\mathbb{R}$: $\ker(D) = \mathbb{R}$;
  \item $H$-field: if $x > \mathbb{R}$, then $D(x) > 0$.
  \end{enumerate}
\end{defn}

Recall that surreal numbers can be written in normal form as
$$
x = \sum_{i<\lambda}\omega^{x_i} \cdot r_i =
\sum_{i<\lambda}e^{\gamma_i} \cdot r_i,
$$
with $\gamma_i \in \mathbb{J}$.

A naive attempt to construct a derivation $D$ is to proceed by
induction, using the following identity, which is a consequence of
rules (2), (3) and (4):
\begin{equation}
  \label{eq:inductive-D}
  D(x) = D\left(\sum_{i<\lambda}e^{\gamma_i} \cdot r_i\right) =
  \sum_{i<\lambda}e^{\gamma_i} \cdot r_i \cdot D(\gamma_i).
\end{equation}

There are two obstructions, however, to such approach: the existence of log-atomic numbers, and finding a way of working by induction.

\subsection{Log-atomic numbers}
\label{sec:log-atomic-numbers}

\begin{defn}
  A surreal number $x \in \no$ is \textbf{log-atomic} if for all
  $n \in \mathbb{N}$ there exists $y_n \in \no$ such that
  $$
  \underbrace{\log(\dots(\log}_{n\textrm{ times}}(x))\dots) =
  \omega^{y_n}.
  $$
  We let $\mathbb{L}$ be the class of all log-atomic numbers.
\end{defn}

For instance, the ordinal numbers $\omega$ and $\epsilon_0$ are
log-atomic numbers; in fact, all $\epsilon$-numbers are log-atomic
\cite{gonshor_surreal}. If $x$ is log-atomic, then
\eqref{eq:inductive-D} does not provide information on its own
regarding the value of $D(x)$; for instance, it does not tell anything
about the value we should give to $D(\epsilon_0)$. It is useful at
this point to classify log-atomic numbers. The first step in this
direction was done in \cite{matu-kuhlm:surreel-transseries-EL} using a
strategy inspired by Conway's definition of the $\Omega$-map.

\begin{defn}[{\cite{matu-kuhlm:surreel-transseries-EL}}]
  Given two positive infinite surreal numbers
  $x, y \in \no^{>\mathbb{N}}$, we say that $x$ and $y$ have the same
  \textbf{exp-log class}, and we write $x \sim^K y$, if there exists
  some $n \in \mathbb{N}$ such that
  $$
  \underbrace{\log(\dots(\log}_{n\textrm{ times}}(x))\dots) < y <
  \underbrace{\exp(\dots(\exp}_{n\textrm{ times}}(x))\dots).
  $$
\end{defn}

A \textbf{$\kappa$-number} is a surreal number that is the simplest in
its own exp-log class. Just as monomials can be parametrised by the
$\Omega$-map, $\kappa$-numbers can be similarly parametrised by the
\textbf{$\kappa$-map}. For the sake of notation, for
$n \in \mathbb{N}$, we define inductively $\log_0(x) = x$,
$\log_{n+1}(x) = \log(\log_n(x))$, and similarly $\exp_0(x) = x$,
$\exp_{n+1}(x) = \exp(\exp_n(x))$. Define
$$
\kappa_x := \{n, \exp_n(x') \mid \log_n(x'') \}
$$
as $n$ varies in $\mathbb{N}$.

We have $\kappa_0 = \omega$, $\kappa_1 = \epsilon_0$. Each
$\kappa$-number is of the form $\kappa_x$ for some $x \in \no$. In
\cite{matu-kuhlm:surreel-transseries-EL} the authors gave a detailed
description of the sign sequence of $\kappa_x$ in term of the sign
sequence of $x$, which in turn can be used to show that all
$\kappa$-numbers are log-atomic.

\begin{thm}[\cite{matu-kuhlm:surreel-transseries-EL}]
  All $\kappa$-numbers, and their images via $\exp_n$, $\log_n$, are
  log-atomic.
\end{thm}

However, the converse does not hold. To capture all log-atomic
numbers, one needs to use the finer notion of \emph{level} adapted
from Hardy fields \cite{rosenlicht:growth,
  marker-miller:levelled-o-min}.

\begin{defn}
  Given two positive infinite surreal numbers
  $x, y \in \no^{>\mathbb{N}}$, we say that $x$ and $y$ have the same
  \textbf{level}, and we write $x \sim^L y$, if there exists some
  $n \in \mathbb{N}$ such that
  $$
    \log_n(x) \sim \log_n(y).
  $$
\end{defn}

As for exp-log classes, a \textbf{$\lambda$-number} is a surreal
number that is the simplest in its own level. Again, $\lambda$-numbers
can be parametrised via the \textbf{$\lambda$-map}:
$$
\lambda_x := \left\{n, \exp_n(n \cdot \log_n(x')) \mid
  \exp_n\left(\frac{1}{n} \cdot \log_n(x'')\right) \right\}
$$
as $n$ varies in $\mathbb{N}^{*}$.

We have $\lambda_0 = \omega$, $\lambda_1 = \exp(\omega)$,
$\lambda_{\omega} = \epsilon_0$. Each $\lambda$-number is of the form
$\lambda_x$ for some $x \in \no$.

\begin{thm}[\cite{ber-man:surreal-derivations}]
  All $\lambda$-numbers are log-atomic, and all log-atomic numbers are
  $\lambda$-numbers.
\end{thm}

Moreover, the identity $\lambda_{x + 1} = \exp(\lambda_x)$ holds for
all $x \in \no$ \cite{adh:surreal-universal}, and comparing this
equation with the definition of $\kappa$-number, one can verify that
$\lambda_x$ is a $\kappa$-number if and only if $x \in \mathbb{J}$.

Note that a priori, surreal derivations can take almost arbitrary
values on $\mathbb{L}$, subject to very few restrictions: for all
$\lambda, \mu \in \mathbb{L}$, one must have
$D(\exp(\lambda)) = \exp(\lambda)D(\lambda)$, and if $\lambda < \mu$,
then $0 < D(\lambda) < D(\mu)$. A further restriction imposed by the
definition of surreal derivation is that one must have
$\log(D(\lambda)) - \log(D(\mu)) \prec \max\{\lambda,\mu\}$ for all
distinct $\lambda, \mu \in \mathbb{L}$.

As in \cite{ber-man:surreal-derivations}, we now define
$\partial_{\mathbb{L}} : \mathbb{L} \to \no$ as the \emph{simplest}
function respecting the above restrictions. It turns out that such
function $\partial_{\mathbb{L}}$ can be calculated quite
explicitly. Given a surreal number $x \in \no$, let $\alpha \in \on$
be the least purely infinite ordinal number such that
$\lambda_{-\alpha} \geq -x$; then
\begin{multline*}
  \partial_{\mathbb{L}}(\lambda_x) = \exp\left(\sum_{i =
      1}^{\infty}\log_i(\lambda_x) - \sum_{\beta < \alpha + 1}
    \lambda_{-\beta}\right) = \\
  = \exp\left(\sum_{i = 1}^{\infty}\log_i(\lambda_x) - \sum_{\beta <
      \gamma + 1} \sum_{i = 1}^{\infty}\log_i(\kappa_{-\beta})\right)
\end{multline*}
where $\gamma$ is the unique ordinal such that
$\lambda_{-\alpha} = \kappa_{-\gamma}$. Note for instance that we
obtain
$\partial_{\mathbb{L}}(\omega) = \partial_{\mathbb{L}}(\lambda_0) =
1$.

\subsection{Extending $\partial_\mathbb{L}$ to $\no$}
\label{sec:extend-part-no}

Let $\mathbb{R}\langle\langle\mathbb{L}\rangle\rangle$ be the smallest
subfield of $\no$ containing $\mathbb{R}$, $\mathbb{L}$, closed under
$\exp$, $\log$, and infinite sums. It is immediate from the definition
that if $\partial_{\mathbb{L}}$ extends to a surreal derivation, then
its values on $\mathbb{R}\langle\langle\mathbb{L}\rangle\rangle$ are
already determined by the values of $\partial_{\mathbb{L}}$ on
$\mathbb{L}$. Since $\mathbb{R}\langle\langle\mathbb{L}\rangle\rangle$
can be constructed inductively from $\mathbb{R}$ and $\mathbb{L}$ by
taking successive closures, \eqref{eq:inductive-D} does provide an
inductive definition of such values. One still has to check that the
definition is well posed, namely that the infinite sums appearing on
the right hand side are indeed summable. This can be done, for
instance, by verifying that $\partial_{\mathbb{L}}$ can be extended
first to $\mathbb{R}((\mathbb{L}))$, and then by applying the results
of \cite{schm01}.

The fact that \eqref{eq:inductive-D} provides an inductive definition
of the extension of $\partial_{\mathbb{L}}$ to
$\mathbb{R}\langle\langle\mathbb{L}\rangle\rangle$ is due to the fact
that $\mathbb{R}\langle\langle\mathbb{L}\rangle\rangle$ satisfies the
condition (ELT4) as mentioned in Section \ref{sec:transseries}; in
fact, $\mathbb{R}\langle\langle\mathbb{L}\rangle\rangle$ is the
largest subfield of $\no$ that is also a field of EL-series in the
sense of \cite{matu-kuhlm:surreel-transseries-EL} (see
\cite{ber-man:surreal-derivations}).

However, $\no$ is strictly larger than
$\mathbb{R}\langle\langle\mathbb{L}\rangle\rangle$, as it contains
numbers as
$$
x = \exp(\omega + \exp(\log_2(\omega) + \exp(\log_4(\omega) + \dots)))
$$
which are solutions to the fixed point equation
$f = \exp(\omega + \log_{2}(f))$ and which do not verify (ELT4).

On the other hand, $\no$ does satisfy the weaker combinatorial
principle (T4) isolated in \cite{schm01} and is therefore a field of
transseries (see Definition \ref{sec:transseries}). Applying
essentially the same technique that Schmeling used to treat
exponential extensions, one can prove that Condition (T4) guarantees
that the right hand side of \eqref{eq:inductive-D} remains summable
\cite{ber-man:surreal-derivations}, so $\partial_{\mathbb{L}}$ does
extend to a surreal derivation. While the extension of
$\partial_{\mathbb{L}}$ from $\mathbb{L}$ to
$\mathbb{R}\langle\langle\mathbb{L}\rangle\rangle$ is unique, its
further extension to $\no$ may not be. Still, there is a ``simplest''
extension $\partial : \no \to \no$.

\begin{thm}[\cite{ber-man:surreal-derivations}]
  There exist several surreal derivations, among which a ``simplest''
  one $\partial : \no \to \no$ extending $\partial_{\mathbb{L}}$.
\end{thm}

\subsection{Universality of $(\no, \partial)$}
\label{sec:theory-no-partial}

It turns out that the derivation $\partial$ has very good
properties. Since $\partial$ is a surreal derivation with
$\partial(\omega) = 1$, \emph{$(\no, \partial)$ is an $H$-field with
  small derivation}.

We now note that the asymptotic couple of the field $(\no, \partial)$
is actually $(\mathbb{J}, \partial)$, up to reversing the ordering of
$\mathbb{J}$, and that $LD(\mathbb{J}^{\neq 0})$ has no infimum in
$\mathbb{J}$ (recall that $\mathbb{J}$ is the non-unital ring of
purely infinite numbers). For the latter, it suffices to check that
$LD(\mathbb{L}) = \partial(\mathbb{L})
= \partial_{\mathbb{L}}(\mathbb{L})$ has no infimum in
$\mathbb{J}$. By the aforementioned Theorem~\ref{thm:asymp-integr} by
Rosenlicht, we obtain the following:

\begin{prop}
  The differential field $(\no, \partial)$ admits asymptotic
  integration.
\end{prop}

On the other hand, since $\no$ can be presented as a union of
spherically complete fields, each one closed under asymptotic
integration, one can deduce that $(\no, \partial)$ is closed under
integration using F.-V. Kuhlmann's
Theorem~\ref{thm:integr-spher-complete}. Since $\no$ also has an
exponential function that by construction of $\partial$ is compatible
with the derivation, this implies the following:

\begin{thm}[\cite{ber-man:surreal-derivations}]
  The differential field $(\no, \partial)$ is Liouville-closed.
\end{thm}

A similar, although way more subtle argument by Aschenbrenner, van den
Dries and van der Hoeven has been used to show that $(\no, \partial)$
is in fact a model of the theory of LE-series, and in turn that it is
universal among all H-fields with constant field $\mathbb{R}$ and
small derivation \cite{adh:surreal-universal}.

\begin{thm}[\cite{adh:surreal-universal}]
  The field $(\no, \partial)$ is an elementary extension of the
  $H$-field of $LE$-series (when identifying the $LE$-series $x$ with
  the number $\omega$).
\end{thm}

The strategy to prove that $(\no, \partial)$ is an elementary
extension of LE-series is in fact the same the authors of
\cite{adh:mod-theo-transseries} use to prove that the field of
LE-series is a model of the theory of
Theorem~\ref{thm:model-complete}. Say that an $H$-field $(K, D)$, with
asymptotic couple $(\Gamma, LD)$, is \textbf{grounded} if
$LD(\Gamma^{\neq 0})$ has a maximum. If one can present an $H$-field
closed under integration as a directed union of spherically complete,
grounded sub-$H$-field, then the $H$-field is in fact both
$\upomega$-free and Newtonian \cite[Cor.\ 11.7.15,
Thm. 15.0.1]{adh:mod-theo-transseries}. The first result follows from
the fact that groundedness implies rather directly
$\upomega$-freeness, which in turn is preserved by directed
unions. Newtonianity follows by a more delicate argument: if the
derivation is surjective in a tower of spherically complete grounded
$H$-fields, then for any asymptotic differential equation defined at
some stage of the tower, a solution in the valuation ring appears few
steps later in the tower (in fact, under some mild assumptions, the
number of steps needed is the order of the differential equation).

To conclude that $(\no, \partial)$ is a model of the theory of
Theorem~\ref{thm:model-complete}, one writes $\no$ as the union of
$K_{\epsilon}$, for $\epsilon$ running over the ordinal epsilon
numbers (i.e.\ such that $\omega^{\epsilon} = \epsilon$), where
$K_{\epsilon}$ is the set of the surreal numbers of the form
$$
\sum_{i < \lambda} r_i\prod_{j < \mu}
\left(\omega^{\omega^{a_{ij}}}\right)^{s_{ij}}
$$
such that either the sign sequence of each $a_{ij}$ is of length less
than $\epsilon$, or $a_{ij} = -\epsilon$.

Each $K_{\epsilon}$ is a spherically complete field, and one can
verify that it is closed under $\partial$, and it is grounded by
construction. This shows that $(\no, \partial)$ is a model of the
theory of LE-series; by model-completeness of such theory, it is in
fact an elementary extension. Moreover, thanks to the explicit
quantifier elimination described in \cite{adh:mod-theo-transseries},
one can in fact prove that $(\no, \partial)$ is a \emph{universal
  $H$-field}:

\begin{thm}[\cite{adh:surreal-universal}]
  Every $H$-field with small derivation and constant field
  $\mathbb{R}$ -- in particular every Hardy field containing $\R$ --
  can be embedded over $\R$ as an ordered differential field into
  $\no$.
\end{thm}

\bibliographystyle{amsalpha}

\end{document}